\documentclass[11pt,reqno]{amsart}
\usepackage{amsfonts,amsmath,amssymb}
\newtheorem{thm}{Theorem}[section]

\newtheorem{corollary}[thm]{Corollary}

\usepackage{color}
\usepackage{graphicx}
\usepackage{epstopdf}

\title[Universal Hilbert]{%
Universal formula for Hilbert series of minimal nilpotent orbits
}
\author{ A. Matsuo}\address{Graduate School of Mathematical Sciences, University of Tokyo,
Komaba, Tokyo, 153-8914, Japan}
\email{matsuo@ms.u-tokyo.ac.jp}

\author{A.P. Veselov}
\address{Department of Mathematical Sciences,
Loughborough University, Loughborough LE11 3TU, UK  and Moscow State University, Moscow 119899, Russia}
\email{A.P.Veselov@lboro.ac.uk}

\begin{document}

\maketitle

\begin{abstract}

We show that the Hilbert series of the projective variety $X=P(O_{min}),$ corresponding to the minimal nilpotent orbit $\mathcal O_{min},$ is universal in the sense of Vogel: it is written uniformly for all simple Lie algebras in terms of Vogel's parameters $\alpha,\beta,\gamma$ and represents a special case of the generalized hypergeometric function ${}_{4}F_{3}.$
A universal formula for the degree of $X$ is then deduced. 

\end{abstract}


\section{Introduction}

Let $\mathfrak g$ be a complex simple Lie algebra and $\mathcal O_{min}$ be the minimal (in the sense of dimension) non-zero nilpotent orbit in $\mathfrak g \approx \mathfrak g^*.$ It is well-known that such orbit is unique, see \cite{CM}).

Its projective version $X=P(\mathcal O_{min}) \subset P(\mathfrak g)$ is a smooth projective variety, sometimes called {\it adjoint variety} (see e.g. \cite{KOY}). It is the only compact orbit of the corresponding complex group $G$ acting on $P(\mathfrak g).$

These varieties have some remarkable geometric properties. In particular, they can be characterised as compact, simply connected, contact homogeneous varieties \cite{Boothby}, or, under certain assumptions, as Fano contact manifolds \cite{Beau}. Their quantum versions related to Joseph ideals \cite{Joseph} were studied in \cite{LSS}.

In this note we analyse these varieties from the point of view of the {\it universal simple Lie algebra} of Vogel \cite{Vogel}.
This object is still to be defined, but Vogel's parametrisation of simple Lie algebras has already proved to be very useful for representation theory (see e.g. \cite{LM,MSV}).

Recall that the universal Vogel's parameters $\alpha, \beta, \gamma$ can be defined as follows \cite{Vogel}.  Take symmetric part of the tensor square of the adjoint representation of $\mathfrak g$ and decompose it into irreducible representations. Then it turns out that, apart from a trivial representation corresponding to invariant symmetric bilinear form, there will be three of them (two for the exceptional Lie algebras). 

Choose an invariant bilinear form on $\mathfrak g$ and consider the eigenvalues of the corresponding Casimir element on these representations. In Vogel's parametrisation they are $4t-2\alpha, 4t-2\beta, 4t-2\gamma,$  where $t=\alpha+\beta+\gamma,$ which defines the parameters uniquely up to a common multiple (because of the freedom in the choice of the invariant form). If we choose the normalisation in which the negative parameter $\alpha=-2,$ then $t=h^\vee$ is the dual Coxeter number and the corresponding parameters are given in Table 1. 

\begin{table}[h]  
\def\arraystretch{1.3}
\begin{tabular}{|c|c|c|c|c|c|}
\hline
Type & Lie algebra  & $\alpha$ & $\beta$ & $\gamma$  & $t=h^\vee$\\   
\hline    
$A_n$ &  $\mathfrak {sl}_{n+1}$     & $-2$ & 2 & $n+1$ & $n+1$\\
$B_n$ &   $\mathfrak {so}_{2n+1}$    & $-2$ & 4& $2n-3 $ & $2n-1$\\
$C_n$ & $ \mathfrak {sp}_{2n}$    & $-2$ & 1 & $n+2 $ & $n+1$\\
$D_n$ &   $\mathfrak {so}_{2n}$    & $-2$ & 4 & $2n-4$ & $2n-2$\\
$E_6$ &  $\mathfrak {e}_{6}  $    & $-2$ & $ 6$& $ 8$ & $12$\\
$E_7$ & $\mathfrak {e}_{7}  $    & $-2$ & $ 8$& $ 12$ & $18$ \\
$E_8$ & $\mathfrak {e}_{8}  $    & $-2$ & $ 12$& $20$ & $30$\\
$F_4$ & $\mathfrak {f}_{4}  $    & $-2$ & $ 5$& $ 6$ & $9$\\
$G_2$ &  $\mathfrak {g}_{2}  $    & $-2$ & $10/3 $& $8/3$ & $4$ \\
\hline  
\end{tabular}
\vskip0.5ex
\leavevmode
\caption{Vogel's parameters for simple Lie algebras}     
\end{table}
Now a formula for numerical characteristics of a simple Lie algebra is called {\it universal}\/ if it is written in terms of Vogel's parameters.

In this note, we will give some universal formulae related to the minimal nilpotent orbit $\mathcal O_{min} \subset \mathfrak g$ by means of the Hilbert series of its projectivisation $X=P(\mathcal O_{min})\subset P(\mathfrak{g})$. 

Recall that for a projective variety $X \subset \mathbb P^n$ the {\it Hilbert series} $H_X(z)$ is defined as the generating function
\begin{equation}
\label{HSdef}
H_X(z)=\sum_{k=0}^\infty \dim(S(X)_k)z^k,
\end{equation}
where $S(X)=\mathbb C[x_0,\dots, x_n]/I(X)$ is the homogeneous coordinate ring of $X$ and $S(X)_k$ is the component of degree $k.$ 
The dimension $\dim(S(X)_k)$ for large $k$ is written as
$$
\dim(S(X)_k)=h_X(k)
$$
by a polynomial $h_X(x)$ called the {\it Hilbert polynomial}, see \cite{Harris}.

Now, return to the case when $X=P(\mathcal{O}_{min})$ and introduce the parameters given by 
\begin{equation}
\label{param1}
a_1=2b_1+2b_2-3, \,\, a_2=b_1+2b_2-2, \,\, a_3=2b_1+b_2-2,
\end{equation}
\begin{equation}
\label{param2}
b_1=-\frac{\beta}{\alpha},\,\,  b_2=-\frac{\gamma}{\alpha},
\end{equation}
where $\alpha,\beta,\gamma$ are Vogel's parameters given in Table 1. 
Recall that the {\it Pochhammer symbol} $(a)_n$ is defined as
\begin{equation}
\label{poch}
(a)_n=a(a+1)\dots(a+n-1), \quad n=1,2,\dots
\end{equation}
with $(a)_0=1.$
Then our main observation is the following universal formula describing the Hilbert series of $X$. 
\begin{thm} 
\label{mymainthm} 
\sl
The Hilbert series of the projectivisation $X=P(\mathcal O_{min})$ of the minimal nilpotent orbit $\mathcal{O}_{min}$ has the following universal form
\begin{equation}
\label{LMseries}
H_{X}(z)=\sum_{k=0}^{\infty}\left(1+\frac{2k}{a_1}\right)\frac{(a_1)_k (a_2)_k (a_3)_k}{ (b_1)_k (b_2)_k k!}z^k,
\end{equation}
with the parameters $a_1,a_2,a_3,b_1,b_2$ given in terms of Vogel's parameters by (\ref{param1}) and (\ref{param2}). 
\end{thm} 
In the rest of the note we will derive the formula (\ref{LMseries}) and describe its consequences. 
In particular, we will extract universal formulae for the dimension and the degree of the projective variety $X$. 
The explicit values of the parameters $a_1,a_2,a_3,b_1,b_2$ as well as $\dim X$ and $\deg X$ are listed in Table 2. 
\begin{table}[htb]
\arraycolsep=0.4em
\def\arraystretch{1.3}
$
\begin{array}{|c|c|c|c|c|c||c|c|}
\hline
\hbox{Type} &a_1&a_2&a_3&b_1&b_2&
\dim X&\deg X\cr
\hline
A_n&n&n&\frac{n+1}{2}&1&\frac{n+1}{2}&2n-1&\binom{2n}{n}\cr
B_n&2n-2&2n-3&n+\frac{1}{2}&2&n-\frac{3}{2}&4n-5&\frac{4}{2n-1}\binom{4n-4}{2n-2}\cr
C_n&n&n+\frac{1}{2}&\frac{n}{2}&\frac{1}{2}&\frac{n}{2}+1&2n-1&2^{2n-1}\cr
D_n&2n-3&2n-4&n&2&n-2&4n-7&\frac{4}{2n-2}\binom{4n-6}{2n-3}\cr
E_6&11&9&8&3&4&21&151164\cr
E_7&17&14&12&4&6&33&141430680\cr
E_8&29&24&20&6&10&57&126937516885200\cr
F_4&8&\frac{13}{2}&6&\frac{5}{2}&3&15&4992\cr
G_2&3&\frac{7}{3}&\frac{8}{3}&\frac{5}{3}&\frac{4}{3}&5&18\cr
\hline
\end{array}
$
\vskip0.5ex
\leavevmode
\caption{Parameters, dimension and degree of $X=P(\mathcal O_{min}).$}
\end{table}

\section{Universal formula for Hilbert series}

The derivation of the universal formula (\ref{LMseries}) is in fact simply a combination of two important results. 

The first one goes back to Borel and Hirzebruch and sometime is attributed to Kostant, see \cite{Gar, GW}. It claims that as a $\mathfrak g$-module $S(X)$ for $X=P(\mathcal O_{min})$ has a form
\begin{equation}
\label{harmonic}
S(X)=\bigoplus_{k=0}^{\infty}V(k \theta),
\end{equation}
where $\theta$ is the maximal root of $\mathfrak g$  and $V(\lambda)$ is the irreducible representation with the highest weight $\lambda$ (see \cite{Gar}).

The second one is the universal formula for the dimension of $V(k\theta)$ found by Landsberg and Manivel \cite{LM}:
\begin{equation}
\label{LM}
\dim V(k \theta)=
\frac{t-(k-\frac 12)\alpha}{ t+\frac{\alpha}2}
\frac{\binom{-\frac{2t}{\alpha}-2+k }{ k}
\binom {\frac{  \beta-2t }\alpha -1+k }{k }
\binom {\frac{  \gamma-2t }\alpha -1+k }{k } 
}
{
\binom{-\frac {\beta}{\alpha} -1+k }{k }
\binom{-\frac {\gamma}{\alpha} -1+k }{k }
},
\end{equation}
where $\binom{a}{k}$ is a binomial coefficient defined for any $a$ and integer $k\geq 0$ as $$\binom{a}{k }=\frac{(a-k+1)\dots(a-1)a}{k!}.$$

Rewriting (\ref{LM}) in terms of Pochhammer symbols (\ref{poch}) and combining it with (\ref{harmonic}), we have
\begin{equation}
\label{LMpoch}
\dim S(X)_{k}
=\left(1+\frac{2k}{a_1}\right)\frac{(a_1)_k (a_2)_k (a_3)_k}{ (b_1)_k (b_2)_k k!},
\end{equation}
where $a_1,a_2,a_3, b_1,b_2$ are given by (\ref{param1}) and (\ref{param2}).
This gives us formula (\ref{LMseries}). 

Let us now turn to the Hilbert polynomial $h_{X}(x)$ which satisfies
$$
h_{X}(k)=\dim S(X)_{k}
$$
if $k$ is sufficiently large.
To obtain a universal formula for $h_{X}(x)$, apply the formula
$$
(a)_k=\frac{\Gamma(a+k)}{\Gamma(a)},
$$
where $\Gamma(x)$ is the classical Gamma-function \cite{WW}.
Then the right-hand side of (\ref{LMpoch}) becomes
$$
\frac{\Gamma(b_1)\Gamma(b_2)}{\Gamma(a_1)\Gamma(a_2)\Gamma(a_3)}
\left(1+\frac{2k}{a_1}\right)\frac{\Gamma(a_1+k)\Gamma(a_2+k)\Gamma(a_3+k)}{\Gamma(b_1+k)\Gamma(b_2+k)\Gamma(1+k)}.
$$
Since this expression for the parameters corresponding to simple Lie algebras turns out to be a polynomial in $k$, we have
\begin{corollary} 
\label{newmainthm} 
\sl
The Hilbert polynomial of $X=P(\mathcal O_{min})$ has the following universal form
\begin{equation}
\label{HilbPoly}
h_{X}(x)
=\frac{\Gamma(b_1)\Gamma(b_2)}{\Gamma(a_1)\Gamma(a_2)\Gamma(a_3)}
\left(1+\frac{2x}{a_1}\right)\frac{\Gamma(a_1+x)\Gamma(a_2+x)\Gamma(a_3+x)}{\Gamma(b_1+x)\Gamma(b_2+x)\Gamma(1+x)},
\end{equation}
with $h_{X}(k)=\dim(S(X)_{k})$ for all $k\geq 0$. 
\end{corollary} 

The fact that the last formula works for all $k \geq 0$ follows from \cite{GW}. 


\section{Geometric consequences}

Recall that Hilbert polynomial of $X$ determines both the dimension and the degree of $X$ as
\begin{equation}
\label{asym1}
h_{X}(x)=\deg \, X \, \frac{x^{d}}{d!}+ \dots,
\end{equation}
where $d =\dim X$ (see e.g. \cite{GW}).
Applying the limit
$$
\lim_{x\to+\infty}\frac{\Gamma(a+x)}{\Gamma(x)}x^{-a}=1,
$$
to (\ref{HilbPoly}), we have the asymptotic expansion
$$
h_{X}(x)=Ax^{a_1+a_2+a_3-b_1-b_2}+\dots,
$$
where the leading coefficient $A$ is given by 
\begin{equation}
\label{asym}
A=\frac{2\,\Gamma(b_1)\Gamma(b_2)}{a_1\Gamma(a_1)\Gamma(a_2)\Gamma(a_3)}.
\end{equation}
Therefore, we have $\dim X=d=a_1+a_2+a_3-b_1-b_2$ and $\deg X=A\cdot d!$ by (\ref{asym1}). 
After a simple algebra, we obtain 
\begin{corollary} 
\label{thm2} 
\sl 
The dimension and the degree of $X=P(\mathcal O_{min})$ are given by 
\begin{equation}
\label{dimension}
\dim X = 2a_1-1,
\end{equation}
\begin{equation}
\label{deg}
\deg(X)=\frac{2\,\Gamma(2a_1)\,\Gamma(b_1)\Gamma(b_2)}{\Gamma(a_1+1)\,\Gamma(a_2)\Gamma(a_3)}
\end{equation}
respectively. 
\end{corollary} 
Since $2a_1-1=-\frac{4t}{\alpha}-3=2h^\vee-3,$ the dimension is more explicitly given by 
$$
\dim X=2h^\vee-3, 
$$
where $h^\vee$ is the dual Coxeter number of $\mathfrak{g}$.
This recovers the formula due to Wang \cite{Wang}: $$\dim \mathcal O_{min}=2h^\vee-2.$$

Let us compare this with the results of Gross and Wallach \cite{GW}, who derived from the Weyl dimension formula that Hilbert polynomial
$h_{X}(q)$ can be given as a product over positive roots of $\mathfrak g$:
\begin{equation}
\label{GW2}
h_{X}(q)=\prod_{\alpha \in R_+}\left(1+\frac{(\theta, \alpha^\vee)}{(\rho, \alpha^\vee)}q\right),
\end{equation}
where $\rho$ is the half-sum of the positive roots.
The Hilbert series of $X$ has the form
\begin{equation}
\label{GW}
H_{X}(z)=h_{X}\left(z\frac{d}{dz}\right)\frac{1}{1-z}.
\end{equation}
As a result they arrived at the following formula for the degree of $X$:
\begin{equation}
\label{GW3}
\deg(X)=d! \prod_{\alpha}\frac{(\theta, \alpha^\vee)}{(\rho, \alpha^\vee)},
\end{equation}
where the product is taken over positive roots such that $(\theta, \alpha^\vee)\neq 0.$
It would be interesting to deduce from here our formula (\ref{deg}), which is a universal form of the right hand side of (\ref{GW3}).

Let us consider as an example the $A_n$-case, when $X=P(\mathcal O_{min})$ is the hyperplane section of the Segre variety \cite{GW,Harris}. The formula (\ref{deg}) gives in this case
$$
\deg X
= \frac{2\,\Gamma(2n)\Gamma(1)\Gamma(\frac{n+1}{2})}{\Gamma(n+1)\Gamma(n)\Gamma(\frac{n+1}{2})}=\frac{2\cdot (2n-1)!}{n!(n-1)!}=\binom{2n}{n},
$$
which agrees with well-known result (see \cite{Harris}). 

The degrees of other types are presented in Table 2. 
The result in $E_8$-case agrees with the one given in \cite{GW}. 
It is interesting to note that in type $B$ and $D$ the degrees are the corresponding Catalan numbers times 4. 

\section{Hilbert series as hypergeometric function}

Let us show that the Hilbert series of $X=P(\mathcal{O}_{min})$ is in fact a particular case of the generalized hypergeometric function ${}_4F_3$ (cf.\ \cite{Andrews}): 
\begin{equation}
\label{4F3}
{}_4 F_{3}(a_1,a_2,a_3,a_4; b_1, b_2; z)=\sum_{n=0}^{\infty}\frac{(a_1)_n (a_2)_n (a_3)_n (a_4)_n}{(b_1)_n (b_2)_n (b_3)_n}\frac{z^n}{n!}. 
\end{equation}
To see this, let the parameters $a_1,a_2,a_3,b_1,b_2$ be given by (\ref{param1}), (\ref{param2}), and $a_4,b_3$ by
\begin{equation}
\label{param3}
a_4=b_3+1,\,\,b_3=-\frac{2t+\alpha}{2\alpha}=\frac{a_1}{2}.
\end{equation}
Then (\ref{LMpoch}) is rewritten as 
\begin{equation}
\label{LMpochNew}
\dim S(X)_k
=\frac{(a_1)_k (a_2)_k (a_3)_k(a_4)_k}{ (b_1)_k (b_2)_k (b_3)_k k!}, 
\end{equation}
and hence we obtain
\begin{corollary} 
\label{mainthm} 
\sl
The Hilbert series of $X=P(\mathcal O_{min})$ 
is written as
\begin{equation}
\label{universal}
H_{X}(z)={}_4 F_{3}(a_1,a_2,a_3, a_4; b_1, b_2, b_3; z), 
\end{equation}
with the parameters $a_1,a_2,a_3,a_4,b_1,b_2,b_3$ given by (\ref{param1}), (\ref{param2}) and (\ref{param3}). 
\end{corollary} 

Note that the series can also be written in terms of ${}_3 F_{2}$ hypergeometric function as
\begin{equation}
\label{universal1}
H_{X}(z)=\left(1+ \frac{2}{a_1} z \frac{d}{dz}\right) {}_3 F_{2}(a_1,a_2,a_3; b_1, b_2; z).
\end{equation}

We can well use the formula (\ref{universal}) to derive the dimension formula for the variety $X$. 
Recall the differential equation for the hypergeometric function $F={}_4F_3$: 
\begin{equation}
\label{eq:DE}
[D\prod_{j=1}^3(zD+b_j-1)-\prod_{j=1}^4(zD+a_j)]F=0, \quad D=\frac{d}{dz}. 
\end{equation}
The point $z=1$ is a regular singular point (in the sense of analytic theory of differential equations) with the exponents $0, 1, 2$ and
$$
s=b_1+b_2+b_3-a_1-a_2-a_3-a_4
$$
(see e.g. \cite{BH}).
Since the Hilbert series $H_{X}(z), \, X=P(\mathcal O_{min})$, is a rational function with the pole at $z=1$ of order equal to the dimension of $\dim X+1$, we have $\dim X=-s-1=a_1+a_2+a_3+a_4-b_1-b_2-b_3-1=2h^\vee-3$. 

Finally, let us recall that the Hilbert series $H_{X}(z)$ can be written as
$$
H_{X}(z)=\frac{P_X(z)}{(1-z)^{d+1}}
$$
for some palindromic polynomial $P_X(z)$ with $P_X(1)=\deg X$ (see e.g.\ \cite{GW}). 
In our case, the problem of finding such a $P_X(z)$ is a particular case of the connection problem for the differential equation (\ref{eq:DE}) from $z=0$ to $z=1$, but formally we can write it down universally as 
\begin{equation}
\label{universal2}
P_{X}(z)=(1-z)^{d+1}H_{X}(z) \,\, ({\rm mod}\,\, (z-1)^{d+1}),
\end{equation}
where $H_X(z)$ is given by (\ref{universal}).

\section{Concluding remarks}

Our results imply also that the Hilbert series of the ideal $I=I(X)$ determining $X$
$$
H_I(z)=\sum_{k=0}^\infty\dim I_k z^k,
$$
where $I_k$ is the component of $I$ of degree $k,$ is also universal. Indeed we have
\begin{equation}
\label{char}
H_I(z)= \frac{1}{(1-z)^N}-H_X(z),
\end{equation}
and $N=\dim \,\mathfrak g$ is the dimension of $\mathfrak g$, which can be given by Vogel's formula
\begin{equation}
\label{dimenV}
N= \frac{(\alpha-2t)(\beta-2t)(\gamma-2t)}{\alpha\beta\gamma}=\frac{(a_1+2)a_2a_3}{b_1b_2}.
\end{equation}
Since 
$$
\frac{1}{(1-z)^N}=\sum_{k=0}^\infty \binom{-N}{k }(-z)^k=\sum_{k=0}^\infty \frac{(N)_k}{k!} z^k,
$$
this implies the following universal expression for $\dim I_k$:
\begin{equation}
\label{dimenIk}
\dim I_k= \frac{(N)_k}{k!}-\left(1+\frac{2k}{a_1}\right)\frac{(a_1)_k (a_2)_k (a_3)_k}{ (b_1)_k (b_2)_k k!}.
\end{equation}

Note that all our formulae are symmetric in $\beta$ and $\gamma,$ but not in $\alpha$.
A natural question is about the meaning of the corresponding Hilbert series when we permute $\alpha$ with $\beta$ or $\gamma$.
From Landsberg and Manivel \cite{LM} we have an interpretation of the corresponding coefficients as the dimensions of certain virtual $\mathfrak g$-modules, which are zero for large $k$ (see the full list for all types in Section 6 of \cite{LM}). 

According to (\ref{dimension}),(\ref{deg}) the corresponding ``virtual varieties" $Y$ and $Z$ must have degree $0$ and negative dimensions:
$$
\dim Y=-\frac{4t}{\beta}-3, \quad \dim Z=-\frac{4t}{\gamma}-3.
$$
In particular, for $A_n$ type $$\dim Y=-2n-5, \,\, \dim Z=-7.$$ Is there any geometry behind this?

\section{Acknowledgements}
We are grateful to Alexey Bolsinov and Jenya Ferapontov for useful discussions.

\end{document}